\title{A Discretization of Bertrand's Paradox}
\author{Gregory Churchill\thanks{
                  SUNY Oswego}
        \and
        Dakota Williams\thanks{SUNY Oswego}}
        
\documentclass{article}
\usepackage{graphicx}

\let\cal=\mathcal

\def\cC{{\cal C}}

\def\cL{{\cal L}}

\def\bbC{\mathbb{C}}

\usepackage{amssymb,amsmath}
\usepackage{mathrsfs}
\usepackage{latexsym}
\usepackage{amscd}
\usepackage{color}

\begin{document}
\newpage
\maketitle
\begin{abstract}
    In this paper we present a discrete version of Bertrand's paradox, where each sample space is a finite set of chords and a chord is chosen uniformly at random.  

\end{abstract}



\section{Introduction}

Consider an equilateral triangle inscribed in a circle.  What is the probability that a randomly chosen chord is longer than a side of the triangle?  This question is the basis for the famously known Bertrand's paradox, which has captured the curiosity of many over the last 100+ years. In his original work, Bertrand offered three different solutions to his question, all of which seem reasonable to a sophomoric reader. But Bertrand was not attempting to answer an interesting probability problem;  he was making a deeper, broader, philosophical point regarding holes in the probability theory of the time. Nonetheless, many have taken a close look at Bertrand's problem and have attempted to solve it, including not just an argument for their answer but also a philosophical justification for why their perspective is the proper one; see [1 -- 4], for example. In a broad internet search, one finds all sorts of arguments for one solution over another, including various simulations to support their answer.

 This articles seeks not to resolve Bertrand's paradox, but rather to further shroud it in mystery. We offer no solution nor simulation; actually, our work here might argue against the use of such simulations in seeking an answer to Bertrand's problem.  

In this paper we present a discrete version of Bertrand's paradox, initially with the following motivation: In Bertrand's original work, and in much of the literature and discussions around his problem, the experiment/sample space for the random experiment does not consist of chords on the circle. Rather, the sample space consists of points on the circle, or points on the disk, or something else, which then determines chords. In this paper, our random experiments have sample spaces that consist of actual chords (albeit they are finite sample spaces), and a chord is chosen uniformly at random. But we still recreate Bertrand's three paradoxical solutions of 1/2, 1/3, and 1/4 (as limits).  Given the nature of our random experiments and sample spaces (to be discussed more in the Concluding Remarks), we hope the reader feels the same sense of paradox one feels when reading Bertrand's original paradox.  

 Here, we always work with the unit circle $S^1$ centered at the origin of the complex plane.  As such, the legs of an inscribed equilateral triangle each have length $\sqrt{3}$, and we refer to a ``long" chord as one whose length is greater than $\sqrt{3}$ (and a ``short" chord as one whose length is less than or equal $\sqrt{3}$).

\section{Discretization of the Angular Separation Method}

For an integer $n \geq 2$, let $R_n$ denote the set of $n$th roots of unity, that is, $R_n = \{ z \in \bbC : z^n = 1\}$.  For $A, B \in S^1$, we write $\overline{AB}$ to denote the chord of $S^1$ with endpoints at $A$ and $B$.  Recall that the $n$th roots of unity are ``evenly spread" on $S^1$, and if $A, B \in R_n$ are ``adjacent" $n$th roots of unity, then the arc length of the minor arc (which we will call {\it minor arc length} for brevity) in $S^1$ between $A$ and $B$ is $2\pi/n$. Let $\Omega_n$ denote the set of all chords of $S^1$ with both endpoints in $R_n$.  Note that $|\Omega_n| = \binom{n}{2}$.  Fix $P = (1,0)$, and we will consider $\ell_P$, the number of long chords which have an endpoint at $P$.  First observe, for $Q \in R_n$, that $\overline{PQ}$ is long if, and only if, $Q$ lies in the minor arc between $\left(-\frac{1}{2}, \frac{\sqrt{3}}{2}\right)$ and $\left(-\frac{1}{2}, \frac{-\sqrt{3}}{2}\right)$ in $S^1$.

On the one hand, we have 

\begin{equation}
\label{eq:less}
\frac{2\pi}{n}(\ell_P - 1) < \frac{2\pi}{3}
\end{equation}   

\noindent Indeed, let $X \in R_n$ be the point closest to $\left(-\frac{1}{2}, \frac{\sqrt{3}}{2}\right)$  such that $\overline{PX}$ is a long chord and let $Y \in R_n$ be the point closest to  $\left(-\frac{1}{2}, \frac{-\sqrt{3}}{2}\right)$ such that $\overline{PY}$ is a long chord.  (Note:  This means, for example, there is no $A \in R_n$ such that $A$ lies in the minor arc between $X$ and $\left(-\frac{1}{2}, \frac{\sqrt{3}}{2}\right)$)  Consider the minor arc length between $X$ and $Y$, which is less than $2\pi/3$ since $X$ and $Y$ both lie on the minor chord between $\left(-\frac{1}{2}, \frac{\sqrt{3}}{2}\right)$ and $\left(-\frac{1}{2}, \frac{-\sqrt{3}}{2}\right)$. But starting at $X$ and moving counterclockwise toward $Y$, we will count all $\ell_P$ roots of unity that determine long chords with $P$ (including $X$ and $Y$ in our count).  Therefore the minor arc length between $X$ and $Y$ is $(\ell_P - 1)2\pi/n$, and (\ref{eq:less}) follows.

On the other hand, we have 

\begin{equation}
\label{eq:more}
\frac{2\pi}{n}(\ell_P + 1) \geq \frac{2\pi}{3}
\end{equation}   

Indeed, let $W \in R_n$ be the point closest to $\left(-\frac{1}{2}, \frac{\sqrt{3}}{2}\right)$  such that $\overline{PW}$ is a short chord and let $Z \in R_n$ be the point closest to  $\left(-\frac{1}{2}, \frac{-\sqrt{3}}{2}\right)$ such that $\overline{PZ}$ is a short chord.  Then starting at $W$ and moving counterclockwise toward $Z$, we will count $\ell_P + 2$ roots of unity.  The minor arc length between $W$ and $Z$ is $(\ell_P + 1)2\pi/n$, which must be greater than or equal to $2\pi/3$ because the minor arc from $W$ to $Z$ contains the minor arc from $\left(-\frac{1}{2}, \frac{\sqrt{3}}{2}\right)$ to $\left(-\frac{1}{2}, \frac{-\sqrt{3}}{2}\right)$. Hence,  (\ref{eq:more}).      

Thus by (\ref{eq:less}) and (\ref{eq:more}) we have $n/3 - 1 \leq \ell_P < n/3 + 1$, bounds on the number of long chords which have an endpoint at $P$.  But by the symmetry of the roots of unity, it follows that these bounds would apply for any fixed $n$th root of unity; that is, for any fixed $A \in R_n$, if $\ell_A$ denotes the number of long chords which have an endpoint at $A$, then we have $n/3 - 1 \leq \ell_A < n/3 + 1$.   

Observe that the total number $\cal{L}$ of long chords of $S^1$ in $\Omega_n$ is 

\begin{equation}
\cal{L} = \frac{1}{2}\sum_{A \in R_n} \ell_A
\end{equation}

\noindent where we divide by 2 because each long chord was counted twice in the summation.  By (\ref{eq:less}) and (\ref{eq:more}) we have 

\begin{equation}
\label{eq:bounds}
\frac{n}{2}\left(\frac{n}{3} - 1\right)
 \ \leq \ \cal{L} \ < \ \frac{n}{2}\left(\frac{n}{3} + 1\right)
\end{equation}

Finally, let $E$ denote the event of selecting a long chord, where we select a chord uniformly at random from $\Omega_n$.  We have 

\begin{equation}
P[E] = \frac{\cal{L}}{|\Omega_n|} = \frac{\cal{L}}{\binom{n}{2}}
\end{equation}

\noindent and with (\ref{eq:bounds}) 

\begin{equation}
\frac{\frac{n}{3}-1}{n-1} \leq \ P[E] \ < \frac{\frac{n}{3}+1}{n-1}
\end{equation}

\noindent and therefore $P[E] \rightarrow 1/3$ as $n \rightarrow \infty$.

  \section{Discretization of the Radial Method}

 Once again fix $P = (1,0) \in S^1$.  Let $O$ denote the origin of the complex plane, and let $\overline{OP}$ denote the radius that intersects $P$.  Fix an integer $m \geq 2$ and consider the points 

\begin{equation}
\label{eq:points}
 \left(\frac{1}{m}, 0\right) , \left(\frac{2}{m}, 0\right) , \dots , \left(\frac{m-1}{m}, 0\right) \in \overline{OP}
 \end{equation} 

\noindent We know that each non-origin point in the unit disk is the midpoint for a unique chord of $S^1$.  We also know that a point $(x,0) \in \overline{OP}$ is the midpoint for a long chord if, and only if, $x < 1/2$.  Hence,  $(j/m , 0) \in \overline{OP}$ is the midpoint for a long chord if, and only if, $j < m/2$.  On the other hand, each of the points 

\begin{equation}
\left(\frac{1}{m}, 0\right) , \left(\frac{2}{m}, 0\right) , \dots ,  \left(\frac{ \lceil m/2 \rceil - 1}{m}, 0\right)   
\end{equation}

\medskip

\noindent is the midpoint for a long chord.  Thus, if $\ell_P$ denotes the number of long chords of $S^1$ which have their midpoints on the points of (\ref{eq:points}), we have 

\begin{equation}
\left\lceil \frac{m}{2} \right\rceil - 1 \  \leq \  \ell_P \ < \ \frac{m}{2} 
\end{equation}

\medskip

Of course the bounds above would apply for any fixed $A \in S^1$; that is, if we take $m-1$ evenly spaced points on $\overline{OA}$ like in (\ref{eq:points}), and  
  if $\ell_A$ denotes the number of long chords of $S^1$ which have their midpoints on these points, we have 

\begin{equation}
 \left\lceil \frac{m}{2} \right\rceil - 1 \  \leq \  \ell_A \ < \ \frac{m}{2} 
\end{equation}

\medskip 

Now to define our experiment space:  For $n \geq 2$, let $R_n$ once again denote the set of $n$th roots of unity, and for each $A \in R_n$, take $m-1$ evenly spaced points on $\overline{OA}$, like in (\ref{eq:points}).  Let $\cC_A$ denote the set of chords which have their midpoints on these $m-1$ evenly spaced points on $\overline{OA}$, and let $\Omega_{m,n} = \cup_{A \in R_n}  \cC_A$. Observe that $|\Omega_{m,n}| = (m-1)n$. Select a chord from $\Omega_{m,n}$ uniformly at random, and let $E$ denote the event of selecting a long chord. Let $\cL$ denote the total number of long chords of $S^1$ in $\Omega_{m,n}$ and observe

\begin{equation}
n\left(\left\lceil \frac{m}{2} \right\rceil - 1\right)
 \ \leq \ \cal{L} = \sum_{A \in R_n} \ell_A  \ < \ n\frac{m}{2}
\end{equation}

\noindent Therefore

   \begin{equation}
P[E] = \frac{\cal{L}}{|\Omega_{m,n}|} = \frac{\cal{L}}{(m-1)n}
\end{equation}

\noindent and 

\begin{equation}
\frac{\left\lceil \frac{m}{2} \right\rceil - 1}{m-1} \ \leq \ P[E] \ < \frac{m}{2(m-1)}
\end{equation}

\noindent so $P[E] \rightarrow 1/2$ as $m,n \rightarrow \infty$.

 \section{Discretization of the Within-disk Method}

Fix $m \geq 2$, and consider $m-1$ concentric circles centered at the origin, with radii $1/m, 2/m, \dots , (m-1)/m$. For $1 \leq j \leq m-1$,  we will write $S^{j/m} = \{z \in \bbC : |z| = j/m\}$, the circle centered at the origin with radius $j/m$.  For $n \geq 2$, let $R^{j/m}_n = \{ z \in \bbC : z^n = j/m\} \subseteq S^{j/m}$.  Again we know that each non-origin point in the unit disk is the midpoint for a unique chord of $S^1$, so let $\cC^{j/m}_n$ denote the set of chords whose midpoints are on $R^{j/m}_n$.  Let $\Omega_{m,n} = \cup_{1 \leq j \leq m-1} \cC^{j/m}_n$ and note $|\Omega_{m,n}| = (m-1)n$.  Once again we select a chord from $\Omega_{m,n}$ uniformly at random and let $E$ denote the event of selecting a long chord.  Once again we let $\cL$ denote the total number of long chords of $S^1$ in $\Omega_{m,n}$.  We know that $\cC^{j/m}_n$ is a set of long chords if, and only if, $j/m \leq 1/4$, i.e., $j \leq m/4$.  But as $j$ is an integer, we have $j \leq m/4$ if, and only if, $ j \leq \lfloor m/4 \rfloor$.   Therefore 

$$ \cL \ = \sum_{1 \leq j \leq \lfloor m/4 \rfloor} |\cC^{j/m}_n|  \ = \ n \left\lfloor\frac{m}{4}\right\rfloor $$

\noindent And thus 

$$ P[E] = \frac{\cL}{|\Omega_{m,n}|} \ = \ \frac{n\lfloor m/4 \rfloor}{n(m-1)} \ = \ \frac{\lfloor m/4 \rfloor}{(m-1)}$$

\medskip
\noindent so that $P[E] \rightarrow 1/4$ as $m,n \rightarrow \infty$.

\section{Concluding Remarks}

While we have reached three different probability calculations, the three experiments above have some common elements:  Each sample space is a set of chords, and a chord is chosen uniformly at random (instead of randomly selecting a point or something else that then determines a chord).  The chords in each sample space are ``evenly spaced out" in some sense;  for example, the roots of unity are evenly spaced out on the circle, and the concentric circles in the within-disk method are evenly spaced apart in the unit disk.  Finally, each sample space can be made arbitrarily large.

But since we have three different probability limit values, we can find three different sets of many chords with different proportions of long chords.  But it is worth repeating that these chords were not ``tailored" to give a biased answer; rather, the chords in each sample space were evenly spaced out in some way, providing a sense of ``fairness."  And since a chord was chosen uniformly at random, every chord in the sample space had the same probability of being selected.  To say the least, we hope this paper encourages a healthy skepticism of simulations that testify to an answer to Bertrand's problem.

\section*{About the authors:}
  Dr. Gregory Churchill is an Assistant Professor of Mathematics at SUNY Oswego and Dakota Williams graduated from SUNY Oswego in 2024 with a Mathematics BA.  Dakota completed his Capstone Project with Dr. Churchill as his supervisor, studying Bertrand's paradox.  This project inspired the investigations of this paper.

\subsection*{Gregory Churchill}
   SUNY Oswego,
   Oswego, New York, 13126.
   gregory.churchill@oswego.edu

\subsection*{Dakota Williams}
   SUNY Oswego,
   Oswego, New York, 13126.
   dwilli27@oswego.edu

\end{document}